\documentclass{article}

\usepackage[english]{babel}

\usepackage[letterpaper,top=2cm,bottom=2cm,left=3cm,right=3cm,marginparwidth=1.75cm]{geometry}

\usepackage{amsmath}
\usepackage{graphicx}
\usepackage{cite}
\usepackage[colorlinks=true, allcolors=blue]{hyperref}
\usepackage[noblocks]{authblk}
\usepackage{algorithm}

\usepackage{amsthm}
\usepackage{arydshln}
\newlength\figureheight
\newlength\figurewidth
\usepackage{amssymb}
\usepackage{mathdots}
\usepackage{multirow}
\usepackage{algpseudocode}
\usepackage{float}

\usepackage{tikz}
\usepackage{tikzscale}
\usetikzlibrary{matrix}
\usepackage{pgfplots}
\usepackage{pdfpages}
\usepackage{hyperref}
\usepackage{mathtools}
\usepackage[noblocks]{authblk}
\usepackage{verbatim}
\usepackage{tikz}
\usetikzlibrary{plotmarks}
\usepackage{subcaption}

\newtheorem{theorem}{Theorem}[section]

\newtheorem{definition}{Definition}[section]
\newtheorem{problem}{Problem}[section]
\newtheorem{proposition}{Proposition}[section]

\newtheorem{example}{Example}[section]

\title{Constructing Sobolev orthonormal rational functions via an updating procedure}
\author[1a*]{Amin Faghih}
	\author[1b]{Marc Van Barel}
	\author[2c]{Niel Van Buggenhout}
	\author[1d]{Raf Vandebril}
 \affil[1]{Department of Computer Science, KU Leuven, Leuven, Belgium}
	
	\affil[2]{Department of Mathematics, Universidad Carlos III de Madrid, Madrid, Spain}

\affil[a*]{Corresponding author: Amin Faghih \texttt{amin.faghih@kuleuven.be}}
\affil[b]{Marc Van Barel \texttt{marc.vanbarel@kuleuven.be}}
\affil[c]{Niel Van Buggenhout \texttt{nvanbugg@math.uc3m.es}}
\affil[d]{Raf Vandebril \texttt{raf.vandebril@kuleuven.be}}

\begin{document}
\maketitle
\begin{abstract}
In this paper, we generate the recursion coefficients for rational functions with prescribed poles that are orthonormal with respect to a continuous Sobolev inner product.
Using a rational Gauss quadrature rule, the inner product can be discretized, thus allowing a linear algebraic approach.
The presented approach involves reformulating the problem as an inverse eigenvalue problem involving a Hessenberg pencil, where the pencil will contain the recursion coefficients that generate the sequence of Sobolev orthogonal rational functions. This reformulation is based on the connection between Sobolev orthonormal rational functions and the orthonormal bases for rational Krylov subspaces generated by a Jordan-like matrix. An updating procedure, introducing the nodes of the inner product one after the other, is proposed and the performance is examined through some numerical examples.
\end{abstract}
	\renewcommand{\thefootnote}{\fnsymbol{footnote}}

	\section{Introduction}
Orthogonal polynomials are an essential tool in developing and analyzing numerical algorithms \cite{OlSlTo20,LiSt13}.
Classical orthogonal polynomials, e.g., Jacobi polynomials, can be constructed immediately from known expressions \cite{Sz75}.
For polynomials that are orthogonal with respect to an inner product not corresponding to a family of classical orthogonal polynomials, numerical linear algebra techniques have been developed to generate them efficiently and reliably, e.g., \cite{GrHa84}.
Let $p_k\in \mathcal{P}_k = \textrm{span}\{1,t,t^2,\dots, t^k\}$ (the space of polynomials up to degree $k$), then the sequence $\{p_0,p_1,p_2,\dots\}$ is an orthogonal polynomial sequence if the orthogonality conditions $(p_i,p_j) = 0$, for $i\neq j$ and $(p_i,p_i)\neq 0$ are satisfied for a continuous inner product
\begin{equation*}
    ( p,q) = \int_{\Omega} p(t)\overline{q(t)} w(t) dt,
\end{equation*}
with weight function $w(t)$ and $\Omega$ a Jordan curve in the complex plane.
One of the most successful approaches for generating a part of the sequence (up to degree $N-1$) $\{p_0,p_1,\dots, p_{N-1}\}$ starts by discretizing the continuous inner product by, e.g., a quadrature rule with $m$ nodes.
This leads to the discretized inner product
\begin{equation*}
    \langle p,q \rangle_m = \sum_{j=1}^m   \vert w_j\vert^2  p(z_j) \overline{q(z_j)}.
\end{equation*}
If the discretization is exact for polynomials up to degree $2N-2$, i.e., $(t^k,t^\ell) = \langle t^k, t^\ell \rangle_m$ holds for $k+\ell\leq 2N-2$, then the polynomials up to degree $N-1$ orthonormal with respect to $(.,.)$ are identical to those orthonormal with respect to $\langle .,. \rangle_m$.

The discretized inner product allows for the use of numerical linear algebra techniques. Typically this amounts to solving an inverse eigenvalue problem. Given the weights and nodes of the discretized inner product: construct a structured matrix (e.g., Hessenberg, tridiagonal, unitary Hessenberg, ...) whose eigenvalues are the nodes and whose eigenvectors link to the weights. The outcome is the structured matrix whose elements are the recursion coefficients for generating the orthogonal polynomials. There are two classical manners of solving the problem: Krylov subspace methods and updating procedures.

For non-classical orthogonal polynomials on the real line, their construction is best done by the updating procedure of Gragg and Harrod \cite{GrHa84} or the continued fraction procedure of Laurie \cite{La99}. A method based on Krylov subspaces is also proposed, namely the Lanczos iteration \cite{dBGo78}. For orthogonal polynomials on the unit circle, an effective method is proposed by Reichel, Ammar, and Gragg \cite{ReAmGr91}.
Apart from generating orthogonal polynomials, these techniques are very effective for solving least squares problems \cite{Re91,BrNaTr21, FaVBVBVa24, VBVBVa23}, where a polynomial approximation to a given set of function evaluations is sought. When additional information underlying the least squares problem is known, this can be taken into account.

If the additional information includes the values of the derivatives at some points of the sought-after function, Sobolev orthogonal polynomials are more appropriate.
Sobolev orthogonal polynomials (SOPs) form a sequence of polynomials orthogonal with respect to a Sobolev inner product.
A \emph{(diagonal) Sobolev inner product} is defined on the space of polynomials $\mathcal{P}$ using $s+1$ finite positive Borel measures $\{\mu_r\}_{r=0}^s$ as:
		\begin{equation}\label{eq:SobInprod}
			(p,q)_S = \sum_{r=0}^{s} \int_{\Omega} p^{(r)}(t) \overline{q^{(r)}(t)} d \mu_r(t),
		\end{equation}
where the support of $\mu_r$, denoted by $\textrm{supp}(\mu_r),$ is a compact subset of the complex plane $\mathbb{C}$.
The study of SOPs is an active field of research \cite{MaXu15,MaZu24,FeMaPePi23,HSMaMaSR23,YuWaLi19}, where the numerical generation of SOPs can be used to study the behavior of the zeros of SOPs and as a tool for solving differential equations.
Methods for generating SOPs are proposed by Gautschi and Zhang \cite{GaZh95}, and for sequentially dominated diagonal Sobolev inner products by Van Buggenhout \cite{VB23}.
Sequentially dominated Sobolev inner products satisfy $\textrm{supp}(\mu_r)\subset \textrm{supp}(\mu_{r-1})$ and $d\mu_r = f_{r-1} d_{\mu_{r-1}}$, with $f_{r-1}\in L_{\infty}(\mu_{r-1})$, $r=1,\dots, s$.
Specialized methods for specific families of SOPs exist and usually outperform these general methods \cite{LaMaVD25}.

If the additional information is the location of a singularity or branch cut, then a basis of rational functions is appropriate to solve the least squares problem \cite{VBVBVa23,BuVBVG04,BeGu15,BeGu17}.
We assume that the poles $\{\xi_k\}_{k\geq 1}$ which capture the singularity or branch cut are known.
Then we include the rational functions with these poles in our space of basis functions.
A sequence $\{r_0,r_1,r_2,\dots\}$ of rational functions with poles $\Xi = \{\xi_k\}_{k\geq 1}$, of the form
\begin{equation*}
	r_n(t)=\dfrac{\mathcal{P}_{n}(t)}{q(t;\Xi)}=\mathcal{R}_n^{\Xi}, \quad \text{with } q(t;\Xi)=\prod_{k=1}^{n}(t-\xi_{k}),
\end{equation*}
is called a sequence of orthogonal rational functions (ORFs) if it satisfies the orthogonality conditions $(r_i,r_j) = 0$, for $i\neq j$ and $(r_i,r_i)\neq 0$. A finite sequence of rational functions orthogonal with respect to $(.,.)$ can be constructed via Krylov subspace methods or updating procedures \cite{VBVBVa22,VBVBVa23}.

If the additional information is both information about singularities or branch cuts and information in the form of point values at derivatives, then we require a basis of rational functions that is orthogonal with respect to a Sobolev inner product, i.e., a basis of Sobolev orthogonal rational functions (SORFs).
In a previous work \cite{FaVBVBVa24} we proposed a Krylov subspace method based on the rational Arnoldi procedure for generating this basis.
In this paper we propose an updating procedure which is more memory efficient than the rational Arnoldi procedure since it does not require to store the whole basis. Furthermore the updating procedure is more flexible as it allows to reuse the existing recurrence relations to build new recurrences under modest modifications, such as adding or deleting nodes, of the inner product. The latter is not possible for the Krylov approach, which would require a costly restart from scratch.

Table \ref{table:methods} provides references for the different methods for all the functions discussed above. The algorithm proposed in this paper completes the table. Note that we have omitted moment-based algorithms \cite{JaRe93,Ga04book}, and biorthogonal functions \cite{VBVBVa22,VB21} from this table.

\begin{table}[t]
\begin{tabular}{l|lll}
                              & Krylov & Updating & Remarks                                                          \\ \cline{1-4}
OPs real line          &  \cite{dBGo78,Re91, Ga82,Ga85}            &       \cite{GrHa84,Re91}              & Krylov procedure, i.e., Lanczos iteration \cite{La50}
\\
OPs unit circle        & \cite{ReAmGr91}                   &                \cite{AmGrRe91b, AmGaRe92}     &                                                                   \\
OPs complex plane          &      \cite{BrNaTr21,VB21}      &               \cite{VB21}      &                 Krylov procedure, i.e., Arnoldi iteration \cite{Ar51}                                                  \\
Sobolev OPs                   &    \cite{VB23}                  &          \cite{VB23}           & Sequentially dominated Sobolev inner products   \\
ORFs &            \cite{VBVBVa22}          &     \cite{VBVBVa22}                &
\\
Sobolev ORFs                  &    \cite{FaVBVBVa24}      &    This paper               &                             Sequentially dominated Sobolev inner products
\end{tabular}
\caption{The summary of the numerical procedures for generating orthogonal polynomials (OPs), Sobolev orthogonal polynomials (SOPs), orthogonal rational functions (ORFs) and Sobolev orthogonal rational functions (SORFs). We list only procedures based on Krylov subspaces and updating.
}
\label{table:methods}
\end{table}

The paper is arranged as follows. Section \ref{sec:theory} is dedicated to reformulating the problem of generating the SORFs to an equivalent matrix problem. In Section \ref{sec:algorithm} we propose the updating procedure to solve the inverse problem.
The numerical experiments in Section \ref{sec:numerics} validate our proposed algorithm and we compare its performance to the rational Arnoldi iteration  \cite{FaVBVBVa24}.
    \section{Inverse eigenvalue problem formulation}
    \label{sec:theory}
    The problem of generating a finite sequence of Sobolev orthogonal rational functions can be equivalently stated as an inverse eigenvalue problem.
    There is some freedom in choosing normalization of these functions, in this paper we consider Sobolev orthonormal rational functions.
    We now state the main problem formally.
    \begin{problem}[Generating Sobolev orthonormal rational functions (SORFs]\label{prob:wHrLS}
	Given a sequentially dominated diagonal Sobolev inner product $(.,.) _S$, as in \eqref{eq:SobInprod}, and a set of poles $\Xi = \{\xi_k\}_{k=1}^{N-1}$, satisfying ${\rm supp}(\mu_r)\cap \Xi = \emptyset$ for all measures $\mu_r$ that appear in $(.,.)_S$. Compute a sequence of rational functions $\{r_k\}_{k=0}^{N-1} $, $r_k \in \mathcal{R}_{k}^{\Xi}$ where
	\begin{equation*}
		( r_k,r_h)_{\rm S}= \left\{\begin{Array}{	cc}
			0, \quad k \ne h,\\
			1, \quad k=h.
		\end{Array}\right.
    \end{equation*}
	\end{problem}
    We use an example to illustrate how to translate this into a discrete problem, which we will solve.
    In this example, we use a Gegenbauer-Sobolev inner product, which has been studied before in the context of Sobolev orthogonal polynomials \cite{MaPePi94} and is an example of a Sobolev inner product that consists of a symmetrically coherent pair \cite{IsKoNoSS91}.
    \begin{example}
        The Gegenbauer-Sobolev inner product with parameters $\mu>-1$ and $\lambda\geq 0$, is defined as follows
        \begin{equation}\label{eq:inprod_example}
            ( p,q)_{{\rm S},\mu} = \int_{-1}^1 p(t)q(t) (1-t^2)^\mu dt + \lambda \int_{-1}^1 p'(t)q'(t) (1-t^2)^\mu dt .
        \end{equation}
    \end{example}

    In order to discretize a Sobolev inner product for the generation of orthogonal rational functions, we require a \textit{rational Gauss quadrature rule} \cite{Ga04book,DeBu12,LoReWu08}.
    This is a rule that is exact on a space of rational functions (as opposed to classical Gauss qaudrature rules, which are exact on a space of polynomials).
    Given an inner product and a set of poles, Gautschi \cite{Ga04book} provided conditions under which a rational Gauss rule exists, as well as an algorithm to compute this rule.
    For a sequentially dominated diagonal Sobolev inner product $(.,.)_S$ as in \eqref{eq:SobInprod}, applying a $\sigma$-node rational Gauss quadrature rule results in a discretized Sobolev inner product of the form
        \begin{equation}\label{eq:DiscretizedInProd}
			\langle r_k,r_h\rangle _S = \sum_{j=1}^\sigma  \sum_{i=0}^{s_{j}} \vert w_j\vert^2 \left\vert \frac{\prod_{r=1}^{i}\alpha_{r}^{(j)}}{i!} \right\vert^2 r^{(i)}_{k}(z_j) \overline{r^{(i)}_{h}(z_j)},
		\end{equation}
    where $s_{j}\leq s$ is the derivative corresponding to $z_{j}$, and the $\alpha_{r}^{(j)}$ are non-zero and allow to vary the weights on the derivatives. The weights are positive if the related poles are simple and real or simple complex conjugates.
    For a proof and more general conditions on the existence of rational Gauss rules with positive weights, we refer to \cite[pp.182-183]{Ga04book} and \cite{DeBu12} and references therein.

    In order to generate a sequence of SORFs for \eqref{eq:SobInprod} with poles $\Xi = \{\xi_k\}_{k=1}^{N-1}$, we require a discretized Sobolev inner product \eqref{eq:DiscretizedInProd} that is exact on the space $\mathcal{R}_{(s+1)(N-1)}^{\Xi^{s+1}}\times \mathcal{R}_{(s+1)(N-1)}^{\Xi^{s+1}}$, where we use the notation $\Xi^{j} = \{ \underbrace{\Xi,...,\Xi}_{j \text{ times}}\}$.
    A rational Gauss quadrature rule with $2(s+1)(N-1)$ finite poles, $\sigma = (s+1)(N-1)+1$ nodes and weights suffices to have exactness on $\mathcal{R}_{(s+1)(N-1)}^{\Xi^{s+1}}\times \mathcal{R}_{(s+1)(N-1)}^{\Xi^{s+1}}$.
    This exactness guarantees that $\langle f,g\rangle_S = (f,g)_S$ for $f,g\in\mathcal{R}_{(s+1)(N-1)}^{\Xi^{s+1}}$. This allows the reformulation of the continuous problem as a discrete, algebraic problem.

    \begin{example}
    \label{ex:gegenbauer_poles}
        Using the Gegenbauer-Sobolev inner product \eqref{eq:inprod_example}, we consider the $2M$ poles $\Xi = \{\xi_1,\xi_2,\dots,\xi_{2M}\} = \{-\omega,\omega,-2\omega,2\omega,\dots, -M\omega,M\omega \}$, with $M \in \mathbb{N}$, and $\omega>1$.
        Any set of poles can be used, as long as it does not overlap with the support of the Gegenbauer measure, i.e., the interval $[-1,1]$. Let us denote the $\sigma= 4M+1$ nodes and weights of this rule as $\{z_{j,\mu},w_{j,\mu}\}_{j=1}^\sigma$.
        These can be computed as described by Gautschi \cite[Example 3.31]{Ga04book}, or by using the Chebfun package \cite{DrHaTr14} for MATLAB.
        Discretization of \eqref{eq:inprod_example} using the quadrature rule for both integrals results in
        \begin{equation}\label{eq:discretizedGegenbauerSobolevInprod}
            \langle p,q\rangle_{S,\mu} = \sum_{j=1}^\sigma  \vert w_{j,\mu}\vert^2 p(z_{j,\mu} ) q(z_{j,\mu}) + \lambda \sum_{j=1}^\sigma  \vert w_{j,\mu}\vert^2 p'(z_{j,\mu} ) q'(z_{j,\mu}).
        \end{equation}
    \end{example}

    Next we relate the space of rational functions $\mathcal{R}_{n-1}^{\Xi}$ to rational Krylov subspaces \cite{Ru84} in order to reformulate Problem \ref{prob:wHrLS} as a matrix problem. There are several essentially equivalent forms to define a rational Krylov subspace, we use the following one.
\begin{definition}[Rational Krylov subspace \cite{Ru84}]\label{def:rks}
	Let $A \in \mathbb{C}^{m\times m}$, $v\in \mathbb{C}^m$,  $\Xi = \{\xi_k \}_{k=1}^{n-1}$ with $\xi_k\in\overline{\mathbb{C}}= \mathbb{C} \cup \{ \infty \}$.
	A rational Krylov subspace $\mathcal{K}_{n}(A,v;\Xi)$, $n \leq m$ with poles $\xi_k$ is defined as
\begin{equation*}
\mathcal{K}_{n}(A,v;\Xi)=  \textrm{span}\left\{v, (A-\xi_1 I)^{-1}v,\dots , \prod_{k=1}^{n-1}  (A-\xi_{k} I)^{-1} v \right\}.
\end{equation*}
\end{definition}

A nested orthonormal basis $Q_{n}$ for $\mathcal{K}_{n}(A,v;\Xi)$ satisfies a recurrence relation of the form
\begin{equation}\label{eq:RKS_recurrence}
	AQ_{n} {K}_{n}+k_{n+1,n} q_{n+1}e_{n}^{\top}= Q_{n} {H}_{n}+h_{n+1,n} q_{n+1}e_{n}^{\top},
\end{equation}
where $h_{n+1,n} \geq 0$ (for an appropriate choice of scaling of the vector $q_{n+1}$), $\frac{h_{k+1,k}}{k_{k+1,k}}= \xi_{k} $, $k=1,2,\ldots,n-1$, $({H}_{n},~{K}_{n}) \in\mathbb{C}^{n\times n} \times \mathbb{C}^{n\times n}$ is a Hessenberg pencil, and $q_{n+1} \in \mathbb{C}^{m}$ satisfies the orthogonality condition $q_{n+1}^{H}Q_{n}=\begin{bmatrix} 0 &  \dots & 0 \end{bmatrix}$.
In fact, by construction the pencil is unreduced, meaning that the entries on the subdiagonals of $H_n$ and $K_n$ are never simultaneously zero.

The vectors spanning a rational Krylov subspace, $\textrm{span}\{q_{1}, q_{2},\ldots,q_{n}\}=\mathcal{K}_{n}(A,v;\Xi)$, can be written in terms of rational functions $r_{k} \in \mathcal{R}_{k}^{\Xi}$, i.e., $q_{k} = r_{k-1}(A)v$. The inner product used is the Euclidean inner product
on $\mathbb{C}^{m}$, i.e., for $x,y \in \mathbb{C}^{m}$, $\langle x,y\rangle= y^{H} x$.

For a given discretized rational Sobolev inner product \eqref{eq:DiscretizedInProd}, we can choose a matrix $A = J$ and starting vector $v = w$ such that we get an equivalence between the inner product on the functions and the Euclidean inner product used for defining the orthogonality of the basis vectors in the Krylov subspace.
In the following theorem, we state this equivalence formally.
		\begin{theorem}
        \label{theo:jl}
        \cite[Section 5]{FaVBVBVa24}\label{theorem:KrylovSobOrs}
		Consider the Jordan-like matrix
		\begin{equation}\label{eq:JordanMatrix}
	 J = \left[\begin{array}{cccc}
		J_{s_{1}} &&&\\
		 &J_{s_{2}}&&\\
		 &&\ddots&\\
		&&&J_{s_{\sigma}}\\
	\end{array}\right]\in \mathbb{C}^{m\times m}, \qquad \text{with }J_{s_{j}} = \begin{bmatrix}
			z_{j}&\alpha_{s_{j}}^{(j)}&&\\
			&z_{j}&&\\
			&&\ddots&\alpha_{1}^{(j)}\\
			&&&z_{j}
		\end{bmatrix} \in\mathbb{C}^{(s_j+1)\times (s_j+1)},
	\end{equation}
	and the vector
    \begin{equation}
    \label{eq:weight}
	w=[\;  \overbrace{0 \ldots  0}^{s_{1}} \;  w_{1} \;  \overbrace{0 \ldots  0}^{s_{2}} \; w_2 \; \dots \; \overbrace{0 \ldots  0}^{s_{\sigma}} \;w_\sigma ]^\top \in \mathbb{C}^{m}.
		\end{equation}
Let $ Q_{n}=\begin{bmatrix}
			q_1 & q_2 & \dots & q_{n}
		\end{bmatrix}  \in\mathbb{C}^{m\times n}$ form a nested orthonormal basis for $\mathcal{K}_{n}(J,w;\Xi)$, $n<m$.
		Assume that the rational functions $r_k\in\mathcal{R}^\Xi_k$ are defined on the spectrum of $J$ satisfying $q_k = r_{k-1}(J)w$.
		Then $\{r_k\}_{k=0}^{n-1}$ is the set of Sobolev orthonormal rational functions with respect to the inner product
		\begin{equation*}
			\langle r_k,r_h\rangle _S = \sum_{j=1}^\sigma  \sum_{i=1}^{s_{j}} \vert w_j\vert^2 \left\vert \frac{\prod_{r=1}^{i}\alpha_{r}^{(j)}}{i!} \right\vert^2 r^{(i)}_{k}(z_j) \overline{r^{(i)}_{h}(z_j)}.
		\end{equation*}
	\end{theorem}

    \begin{example}\label{ex:gegenbauer_matrixJ}
        Considering Example~\ref{ex:gegenbauer_poles}, linked to inner product \eqref{eq:discretizedGegenbauerSobolevInprod}, we see that the Jordan blocks $J_{s_{j}}$  and vector $w$ from Theorem~\ref{theo:jl} are
        \begin{equation*}
	J_{s_{j}} = \begin{bmatrix}
			z_{j,\mu}& \sqrt{\lambda}\\
			&z_{j,\mu}\\
		\end{bmatrix} \in \mathbb{C}^{2 \times 2}, \text{ and} \quad w=[\;  0 \;  w_{1,\mu} \; 0 \; w_{2,\mu} \; \dots \; 0 \;w_{\sigma,\mu} ]^\top \in \mathbb{C}^{2\sigma}.
	\end{equation*}
    The poles $\Xi = \{-\omega,\omega,-2\omega,2\omega,\dots, -M\omega,M\omega \}$ of the SORFs influence the quadrature nodes and weights directly through the Gaussian quadrature rule, i.e., to get a correct quadrature rule the weights and nodes are dependent on the poles. The poles themselves must be placed explicitly in the Hessenberg pencil $(H_{n},K_{n})$ in \eqref{eq:RKS_recurrence}, encoded as the ratio of subdiagonal entries. This leads to Problem \ref{IEP} formulated below.
    \end{example}

	The consequence of Theorem \ref{theorem:KrylovSobOrs} is that the recurrence coefficients that generate orthogonal vectors in some rational Krylov subspaces can be directly used to construct Sobolev orthogonal rational functions. These functions will be orthogonal with respect to a discretized inner product that is composed of the eigenvalues of $J$ and elements of $w$, which make up the rational Krylov subspaces. Hence, Problem \ref{prob:wHrLS} can be formulated as a matrix problem using Theorem \ref{theorem:KrylovSobOrs}. It takes the form of an inverse eigenvalue problem (IEP). The considered inverse eigenvalue problem is provided in Problem \ref{IEP}. First let us give the following proposition.
	\begin{proposition}\label{proposition}\cite{FaVBVBVa24,VB23}
		Consider a nested orthonormal basis $Q_{n} \in \mathbb{C}^{m \times n}$ for the rational Krylov subspace $\mathcal{K}_{n}(J,w;{\Psi})$, with $n < m$ and $J \in \mathbb{C}^{m \times m}$ , $w \in \mathbb{C}^{m}$ as in \eqref{eq:JordanMatrix} and \eqref{eq:weight}.
        The set of poles ${\Psi}= \{{\psi}_k\}_{k=1}^{n}$ can be chosen arbitrarily as long as they do not overlap with the spectrum of $J$.
        Let $({H}_{n},~{K}_{n}) \in\mathbb{C}^{n\times n} \times \mathbb{C}^{n\times n}$ denote an unreduced Hessenberg pencil such that $\frac{h_{k+1,k}}{k_{k+1,k}} = {\psi}_k$, and
		\begin{equation*}
			JQ_{n} {K}_{n}+k_{n+1,n} q_{n+1}e_{n}^{\top}= Q_{n} {H}_{n}+h_{n+1,n} q_{n+1}e_{n}^{\top}.
		\end{equation*}
		Then, the recurrence relation for the sequence of SORFs $\{r_{k}\}_{k=0}^{n}$, with $n<m$ and $r_k\in \mathcal{R}_{k}^{{\Psi}}$, orthonormal with respect to $\langle .,.\rangle _S$ \eqref{eq:DiscretizedInProd}, is given by
	\begin{equation}\label{eq3}
		t \begin{bmatrix}
			r_0 & r_1 & \dots & r_{n-1}
		\end{bmatrix} {K}_{n}+k_{n+1,n} r_{n}(t)= \begin{bmatrix}
			r_0 & r_1 & \dots & r_{n-1}
		\end{bmatrix}{H}_{n}+h_{n+1,n} r_{n}(t).
	\end{equation}
		The converse also holds.
	\end{proposition}
	For $n = m$, the recurrence relation \eqref{eq3} terminates, i.e.,
	\begin{equation*}
		JQ_{m} {K}_{m}= Q_{m} {H}_{m},
	\end{equation*}
	with $Q_{m} e_{1} = w/|| w ||_{2}$ and $Q_{m}^H Q_m =I$ . Since $Q_{m}$ is unitary, we have
	$Q_{m}^H J Q_{m} ={H}_{m} {K}_{m}^{-1}$,
	which is the Jordan canonical form of the matrix ${H}_{m} {K}_{m}^{-1}$. Conversely, if the poles $\Psi$,
	the first entries of the eigenvectors $Q_{m} e_1 = w/||w||^{2}$ and the Jordan matrix $J$ are
	known, then the unitary matrix $Q_{m}$ and Hessenberg pencil $({H}_{m},~{K}_{m})$ can be reconstructed.

    Proposition \ref{proposition} provides the essential connection between an orthonormal basis for a rational Krylov space and the orthonormal basis for a rational function space with respect to the discretized inner product $\langle .,.\rangle_S$.
    However, we are interested in an orthonormal basis for rational functions in $\mathcal{R}_{N-1}^{\Xi}$ with respect to the continuous inner product $(.,.)_S$.
    In order to be able to use Proposition \ref{proposition} for the continuous inner product, it is essential that we take $\psi_k=\xi_k$, for $k=1,2,\dots, N-1$, as the ratio of the subdiagonal entries in the Hessenberg pencil, and that the matrix $J$ and vector $w$ are constructed using a rational Gauss quadrature as described above.
	\begin{problem}[Hessenberg pencil IEP]\label{IEP}
		Consider a Jordan-like matrix
			\begin{equation*}
			J = \left[\begin{array}{cccc}
				J_{s_{1}} &&&\\
				&J_{s_{2}}&&\\
				&&\ddots&\\
				&&&J_{s_{\sigma}}\\
			\end{array}\right]\in \mathbb{C}^{m\times m},
		\end{equation*}
		 with blocks $	J_{s_{j}} \in \mathbb{C}^{(s_{j}+1) \times (s_{j}+1)}$,
         with distinct nodes $z_{j} \in \mathbb{C}$, and a vector containing weights $\{w_{j}\}_{j=1}^{\sigma}$
			\begin{equation*}
			w=\begin{bmatrix}      \overbrace{0 \ldots  0}^{s_{1}} & w_{1} &   \overbrace{0 \ldots  0}^{s_{2}} & w_{2} & \dots& \overbrace{0 \ldots  0}^{s_{\sigma}} & w_{\sigma} \end{bmatrix}^\top \in \mathbb{C}^{m}.
		\end{equation*}
	Given a set of poles $\Xi =\{\xi_{k}\}_{k=1}^{N-1}$,
	 $\xi_{k} \in \overline{\mathbb{C}}$. Construct a Hessenberg pencil $({H}_{m},~{K}_{m}) \in\mathbb{C}^{m\times m}$, with
		$\frac{h_{k+1,k}}{k_{k+1}}= \xi_{k} $, $k=1,2,\ldots,N-1$, and $m=\sum_{j=1}^\sigma (s_j+1)$ such that
	\begin{equation}\label{eq:IEP_solution}
		Q_{m}^{H} J Q_{m} {K}_{m}={H}_{m},\quad \text{and} \quad Q_{m} e_{1} = w/|| w ||_{2}.
	\end{equation}
	\end{problem}
	This problem implies, by the rational implicit $Q$ theorem \cite{CaMeVa19_QZ}, that the column span of $Q_{m}$ equals $\mathcal{K}_{m}(J,w;\Psi)$, where $\Psi = \{\psi_k\}_{k=1}^{m-1}$ with $\psi_k=\xi_k$ for $k=1,2,\dots, N-1$.
    In this theorem, the key is to note the correspondence between the poles in the Krylov subspace and the ratio of the subdiagonal elements of the Hessenberg pencil.
    Solving this inverse eigenvalue problem corresponds to computing recurrence coefficients of a sequence of Sobolev orthogonal rational functions in $\mathcal{R}_{m-1}^{\Psi}$, orthogonal with respect to a discretized inner product.
    The first $N$ rational functions form a nested basis for $\mathcal{R}_{N-1}^{\Xi}$ and they are orthonormal with respect to the continuous inner product $(.,.)_S$. Hence, they form a solution to Problem \ref{prob:wHrLS}.

	\section{Updating procedure}\label{sec:algorithm}
	By recasting Problem \ref{prob:wHrLS} into the matrix formulation of Problem \ref{IEP}, we can rely on numerical linear algebra techniques to generate SORFs.
    We propose an updating	procedure by using unitary equivalence transformations with plane rotations.

    The main idea behind the updating procedure is to add nodes one after the other; as each node is linked to a Jordan block, this amounts to adding Jordan blocks one after the other. In fact, we split a discretized Sobolev inner product of interest \eqref{eq:DiscretizedInProd} into the sum of two inner products:
		\begin{eqnarray*}
		\langle r,s\rangle _{\tilde{S}} &= \underbrace{\sum_{j=1}^{\sigma-1}  \sum_{i=1}^{s_{j}} \vert w_j\vert^2 \left\vert \frac{\prod_{r=1}^{i}\alpha_{r}^{(j)}}{i!} \right\vert^2 r^{(i)}(z_j) \overline{s^{(i)}(z_j)}}_{	\langle r,s\rangle _{\hat{S}}}
		&+\underbrace{\sum_{i=1}^{s_{\sigma}} \vert w_\sigma\vert^2 \left\vert \frac{\prod_{r=1}^{i}\alpha_{r}^{(\sigma)}}{i!} \right\vert^2 r^{(i)}(z_\sigma) \overline{s^{(i)}(z_\sigma)}}_{	\langle r,s\rangle _\sigma}.
	\end{eqnarray*}
	This splitting is equivalent to dividing the corresponding matrix and vector into subblocks as follows,
	\begin{equation*}
		\tilde{J} = \left[\begin{array}{ccc|c}
			J_{s_{1}} &&&\\
			&\ddots&&\\
			&&J_{s_{\sigma-1}}&\\ \hline
			&&&J_{s_{\sigma}}\\
		\end{array}\right]=\left[\begin{array}{c|c}
		\hat{J} &\\\hline
		&J_{s_{\sigma}}\\
	\end{array}\right],
	\end{equation*}
	and
		\begin{equation*}
		\tilde{w}= [\;  \overbrace{0 \ldots  0}^{s_{1}} \;  w_{1} \;  \overbrace{0 \ldots  0}^{s_{2}} \; w_2 \; \dots \; \rvert \; \overbrace{0 \ldots  0}^{s_{\sigma}} \;w_\sigma ]^\top
        = [\; \hat{w}\;\rvert\; \overbrace{0 \ldots  0}^{s_{\sigma}} \;w_\sigma ]^\top.
	\end{equation*}
    We assume that the recurrence coefficients for the sequence of SORFs for $\langle r,s\rangle _{\hat{S}}$ and the sequence for $\langle r,s\rangle_\sigma$ are already computed, i.e., the solution to Problem \ref{IEP} is known for both submatrices.
    Starting from these two known partial solutions, we can efficiently construct the Hessenberg pencil for the SORFs for $\langle r,s\rangle _{\tilde{S}}$.
    There are, essentially, four steps, which are summarized here and described in more detail below.
	\begin{enumerate}
		\item Embed the solutions to Problem \ref{IEP} corresponding to the two smaller submatrices of size $m-(s_{\sigma}+1)$ and $s_{\sigma}+1$ into a larger pencil of size $m$. The resulting Hessenberg pencil has the nodes of $\langle r,s\rangle _{\tilde{S}}$ as its eigenvalues.
		
		\item The Hessenberg pencil does not satisfy \eqref{eq:IEP_solution}, since the corresponding basis matrix $Q_{m}$ does not have the correct first column.
        Via multiplication by a unitary matrix we can change the first column such that it is equal to the normalized new weight vector of $\langle r,s\rangle _{\tilde{S}}$.
        As a consequence, \eqref{eq:IEP_solution} will be satisfied, however, the structure of the corresponding recurrence pencil is no longer Hessenberg.
		
		\item  Via unitary equivalence transformations we enforce the Hessenberg structure on the recurrence pencil, while preserving \eqref{eq:IEP_solution}.

        \item The poles $\Psi$ must be introduced on the subdiagonal of the Hessenberg pencil while preserving \eqref{eq:IEP_solution}. This can be done by existing pole swapping techniques \cite{CaMeVa19_QZ}.
\end{enumerate}

\subsection{Embedding of partial solutions}
The recurrence pencils for the SORFs with respect to $\langle .,. \rangle_{\sigma}$ and $\langle .,.\rangle_{\hat{S}}$ are required to start the updating procedure to construct the pencil for the SORFs for $\langle .,. \rangle_{\tilde{S}}$.
These two pencils are then embedded in one large pencil.

For the first inner product $\langle r,s\rangle _{\sigma}$, the Hessenberg pencil IEP is formulated for a single block $J_{s_{\sigma}} \in \mathbb{C}^{(s_{\sigma}+1) \times (s_{\sigma}+1)}$ and corresponding weight vector $w=[\; \overbrace{0 \ldots  0}^{s_{\sigma}} \;w_\sigma ]^\top$.
This Hessenberg pencil can be easily obtained, let $I_{j}$ denote the identity matrix of size $j \times j$, then
\begin{equation*}
	{H}_{s_{\sigma}+1}=J_{s_{\sigma}}^{H}, \quad {K}_{s_{\sigma}+1}=I_{s_{\sigma}+1}, \quad Q_{s_{\sigma}+1}=\begin{bmatrix}
		&&1\\
		&\iddots&\\
		1&&\\
	\end{bmatrix}\in \mathbb{C}^{(s_{\sigma}+1) \times (s_{\sigma}+1)},
\end{equation*}
satisfy the corresponding conditions \eqref{eq:IEP_solution}, i.e., the recurrence relation $J_{s_{\sigma}}Q_{s_{\sigma}+1}{K}_{s_{\sigma}+1}=Q_{s_{\sigma}+1}{H}_{s_{\sigma}+1}$ and the weight condition $Q_{{s_\sigma+1}}e_1 = \frac{w}{\Vert w\Vert_2}$.

For the second inner product $\langle r,s\rangle _{\hat{S}}$, we assume, without loss of generality, that we know the solution, ${\hat{H}},{\hat{K}} \in \mathbb{C}^{m-(s_{\sigma}+1) \times m-(s_{\sigma}+1)}$ and $\hat{Q}\in \mathbb{C}^{m-(s_{\sigma}+1) \times m-(s_{\sigma}+1)}$, to the corresponding Hessenberg pencil IEP.
Without loss of generality, because we can distinguish two cases. First case, if $\hat{J}$ is a single block, the solution is obtained as described above, and we just require to add $s_{1}$ poles to the Hessenberg pencil and swap them to the correct position (The procedure of swapping a pole is presented in Section \ref{sec:addSwapPoles}). In the second case we can build the solution by running the updating procedure described below, starting from a single block and adding one block at a time.

Now that the solutions to Problem \ref{IEP} for both inner products are known, they are embedded in the larger matrices of size $m$,
\begin{equation*}
	\widetilde{Q}=\begin{bmatrix}
		\hat{Q}& \\
		& Q_{s_{\sigma}+1}
	\end{bmatrix}, \quad 	\widetilde{H}=\begin{bmatrix}
	{\hat{H}}& \\
	& J_{s_{\sigma}}^{H}
\end{bmatrix},\quad \text{and} \quad	{\widetilde{K}}=\begin{bmatrix}
{\hat{K}}& \\
& I
\end{bmatrix},
\end{equation*}
where $\widetilde{Q}$ is unitary and the pencil $({\widetilde{H}},~\widetilde{K})$ satisfies the recurrence relation $\widetilde{J}\widetilde{Q}{\widetilde{K}}=\widetilde{Q}{\widetilde{H}}$.
However, the first entries of the eigenvectors are not equal to the given weights $\widetilde{Q}e_{1} \ne  \frac{\tilde{w}}{|| \tilde{w} ||_{2}}$, thus these larger matrices do not satisfy \eqref{eq:IEP_solution}.
The following section describes how this can be achieved by a plane rotation.

\subsection{Adjusting the basis matrix}
A plane rotation is applied to the basis matrix $\widetilde{Q}$, whose action alters the first column such that it equals the normalized weight vector $\frac{\tilde{w}}{|| \tilde{w} ||_{2}}$.
Plane rotations are essentially $2 \times 2$ unitary matrices embedded in a larger identity matrix. They have two parameters $a,b \in \mathbb{C}$ satisfying $\bar{a}a + \bar{b}b = 1$. The plane rotations used for introducing weights into $\widetilde{Q}$ are of the form
\begin{equation}\label{eqqpl1}
	P_{j}=\begin{bmatrix}
		\bar{a}&&-\bar{b}&\\
		&I_{m-j-2}&&\\
		b&&a&\\
		&&&I_{j}
	\end{bmatrix}.
\end{equation}

By taking $a= \frac{|| \hat{w} ||_{2}}{|| \tilde{w} ||_{2}}$ and $b=\frac{w_{\sigma}}{|| \tilde{w} ||_{2}}$ in \eqref{eqqpl1}, the matrix $\widetilde{Q}P^{H}_{s_\sigma}$ satisfies $\widetilde{Q}P^{H}_{s_\sigma}e_{1}=\frac{\tilde{w}}{|| \tilde{w} ||_{2}}$.
In the new basis formed by columns of $\widetilde{Q}P^{H}_{s_\sigma}$, the matrix $\widetilde{J}$ can be represented by the pencil $(P_{s_\sigma}{\widetilde{H}},P_{s_\sigma}{\widetilde{K}})$, this follows from
\begin{equation*}
\widetilde{J}\widetilde{Q}P^{H}_{s_\sigma}P_{s_\sigma}{\widetilde{K}}=\widetilde{Q}P^{H}_{s_\sigma}P_{s_\sigma}{\widetilde{H}}.
\end{equation*}
The resulting pencil $(P_{s_\sigma}{\widetilde{H}},P_{s_\sigma}{\widetilde{K}})$ is no longer a Hessenberg pencil.\\

Using an example, with different orders of derivatives for each block, we will show the structure of the pencil during the updating procedure.
In Example \ref{ex1} we show the structure of the Hessenberg pencil $({\widetilde{H}},{\widetilde{K}})$ and the structure of the pencil $(P_{s_\sigma}{\widetilde{H}},P_{s_\sigma}{\widetilde{K}})$ resulting from correcting the weights in the first column of the basis matrix.
\begin{example}\label{ex1}
	Let us consider the Gegenbauer-Sobolev inner product \eqref{eq:inprod_example} with a given pole $\Xi=\{\xi_{1}\}$. Since $s=1$, the discretized Sobolev inner product \eqref{eq:DiscretizedInProd} is exact for $\sigma=3$. $\{s_{1},s_{2},s_{3}\}=\{1,1,1\}$ is the set of derivatives associated with each node, and from Example \ref{ex:gegenbauer_matrixJ}, we know that $m=6$, thus $\Psi=\{\psi_{k}\}_{k=1}^{5}$ with $\psi_{1}=\xi_{1}$. For adding the block with $s_{3}=1$ to the solution already computed for $s_1,s_2$, the embedded matrix ${\widetilde{H}}$ is mapped onto $P_{1}{\widetilde{H}}$, with the following structure,
	\begin{equation*}
{\widetilde{H}}=	\begin{bmatrix}
		\times & \times & \times &\times&&\\
		\times & \times & \times  &\times&&\\
		 & \times & \times  &\times&&\\
		
		   &  & \times  &\times&&\\
		     &  &   &&\times&\\
		     &  &   &&\times&\times\\
	\end{bmatrix}, \quad
P_{1}{\widetilde{H}}=\begin{bmatrix}
\times & \times & \times & \times & \times &   \\
\times & \times & \times & \times &   &   \\
  & \times & \times & \times &   &   \\
  &   & \times & \times &   &   \\
\times & \times & \times & \times & \times &   \\
  &   &   &   & \times & \times
\end{bmatrix},
	\end{equation*}
	for the second matrix in the pencil, ${\widetilde{K}}$ is mapped onto $P_{1}{\widetilde{K}}$ with the following structure,
		\begin{equation*}
		{\widetilde{K}}=	\begin{bmatrix}
\times & \times & \times & \times &   &   \\
\times & \times & \times & \times &   &   \\
  & \times & \times & \times &   &   \\
  &   & \times & \times &   &   \\
  &   &   &   & \times &   \\
  &   &   &   &   & \times
\end{bmatrix}, \quad
		P_{1}{\widetilde{K}}=\begin{bmatrix}
\times & \times & \times & \times & \times &   \\
\times & \times & \times & \times &   &   \\
  & \times & \times & \times &   &   \\
  &   & \times & \times &   &   \\
\times & \times & \times & \times & \times &   \\
  &   &   &   &   & \times
\end{bmatrix},
	\end{equation*}
where $\times$ denotes a generic non-zero element.
\end{example}
The following subsection describes how the Hessenberg structure can be restored.

\subsection{Restore Hessenberg structure}
Via unitary equivalence transformations the pencil $(P_{s_\sigma}{\widetilde{H}},P_{s_\sigma}{\widetilde{K}})$ is reduced to Hessenberg structure without altering the first column of $\widetilde{Q}P^{H}_{s_\sigma}$.
The following operation is essential, it allows us to eliminate entries in the pencil while preserving the poles that have been introduced (i.e., the ratios of the subdiagonal entries).
\begin{itemize}
\item{\textbf{Operation 1} (Element elimination with pole preservation).}
	
Let us consider the pencil $({\widetilde{H}^{[0]}},{\widetilde{K}^{[0]}})=(P_1{\widetilde{H}},P_1{\widetilde{K}})$ obtained in Example \ref{ex1} and assign the following parameters,
\begin{eqnarray*}
	{\widetilde{H}^{[0]}}=\begin{bmatrix}
\times & \times & \times & \times & \times &   \\
\delta & \times & \times & \times &  0 &   \\
  & \times & \times & \times &   &   \\
  &   & \times & \times &   &   \\
\gamma & \times & \times & \times & \eta &   \\
  &   &   &   & \times & \times
\end{bmatrix} , \quad \quad
	{\widetilde{K}^{[0]}}=  \begin{bmatrix}
\times & \times & \times & \times & \times &   \\
\beta & \times & \times & \times & 0  &   \\
  & \times & \times & \times &   &   \\
  &   & \times & \times &   &   \\
\alpha & \times & \times & \times & \zeta &   \\
  &   &   &   &   & \times
\end{bmatrix}.
\end{eqnarray*}
Our aim is to eliminate the entries $\gamma$ and $\alpha$ in a way that $\psi_{1}=\frac{\delta}{\beta}$ is preserved. Since plane rotations are used for this purpose, we isolate the relevant elements in an equivalent $2 \times 2$ problem, i.e., find parameters $a,b$ and $c, d$ appearing in plane rotations $P$ and $\dot{P}$, respectively, such that
\begin{eqnarray*}
	PA\dot{P}&=&\begin{bmatrix}
		\bar{a} & -\bar{b}\\
		b & a
	\end{bmatrix} \begin{bmatrix}
		\delta& 0\\
		\gamma & \eta
	\end{bmatrix} \begin{bmatrix}
		\bar{c} & -\bar{d}\\
		d & c
	\end{bmatrix}=\begin{bmatrix}
		h & \times\\
		0 & \times
	\end{bmatrix},\\
	PB\dot{P}&=&\begin{bmatrix}
		\bar{a} & -\bar{b}\\
		b & a
	\end{bmatrix} \begin{bmatrix}
		\beta& 0\\
		\alpha & \zeta
	\end{bmatrix} \begin{bmatrix}
		\bar{c} & -\bar{d}\\
		d & c
	\end{bmatrix}=\begin{bmatrix}
		k & \times\\
		0 & \times
	\end{bmatrix},
\end{eqnarray*}
with $\frac{h}{k}=\frac{\delta}{\beta}=\psi_{1}$, and $2\times 2$ blocks $A$ and $B$ extracted from ${\widetilde{H}^{[0]}}$ and ${\widetilde{K}^{[0]}}$, respectively. In this regard, $\dot{P}$ is first constructed in a way that $M \dot{P}e_{1}=\begin{bmatrix} 0& 0\end{bmatrix}^{\top}$ with
\begin{equation*}
	M=\beta A-\delta B=\begin{bmatrix}
		0&0\\
		\times & \times
	\end{bmatrix}.
\end{equation*}
This results in  ${\rm rank}\left(\begin{bmatrix}
    A \dot{P}e_{1} &  B\dot{P}e_{1}
\end{bmatrix}\right)=1$, i.e., the first columns of $A\dot{P}$ and $B \dot{P}$ are colinear. Second, $P$ is selected to make $A \dot{P}$ (or $B \dot{P}$) upper triangular, i.e., $P A\dot{P} e_{1} = \begin{bmatrix} h& 0\end{bmatrix}^{\top}$. Thanks to the colinearity, the same plane rotation $P$ results in $P B \dot{P} e_{1} = \begin{bmatrix} k& 0\end{bmatrix}^{\top}$. For some non-zero constant $u$, we have $h=u\delta$ and $k=u \beta$ and, therefore the ratio of subdiagonal elements $\frac{h}{k}=\frac{u\delta}{u\beta}=\psi_{1}$ is preserved.
Returning to Example \ref{ex1}, the elements on position $(5,1)$ can thus be eliminated and the structure of the resulting pencil is,
\begin{eqnarray}
	\left( \left[\begin{array}{cccccc}
\times & \times & \times & \times & \times &   \\
\times & \times & \times & \times & \times &   \\
  & \times & \times & \times &   &   \\
  &   & \times & \times &   &   \\
  & \times & \times & \times & \times &   \\
\times &   &   &   & \times & \times
	\end{array}\right] ,
    \left[\begin{array}{cccccc}
		 \times & \times & \times & \times & \times &   \\
\times & \times & \times & \times & \times &   \\
  & \times & \times & \times &   &   \\
  &   & \times & \times &   &   \\
  & \times & \times & \times & \times &   \\
  &   &   &   &   & \times
	\end{array}\right]
    \right). \label{eq:struct_Ex3.1}
\end{eqnarray}

Hence, this process can be repeated two more times to eliminate all the entries below the first subdiagonal of the first column.
\end{itemize}

Operation 1 is used column by column, starting at the first column we eliminate all nonzero entries below the subdiagonal, move on to the second column and work up to the column $(m-2)$.
In more detail, let $r$ and $c$ denote the number of the row and column considered, respectively.
For the first column, $c=1$, up to the column $c=(m-1)-(s_{\sigma}+1)$, we eliminate the entries of the rows $r=(m-s_{\sigma})$ to $r=m$. Afterwards, for $c=m-(s_{\sigma}+1),\ldots,m-2$, we eliminate the entries in the rows $r=(c+2),\ldots,m$.
In general, $((m-1)-(s_{\sigma}+1))\cdot (s_{\sigma}+1)+\sum_{j=1}^{\sigma}s_{j}$ entries should be eliminated from $P_{s_\sigma}{\widetilde{H}}$ and $P_{s_\sigma}{\widetilde{K}}$.

The elimination is performed using Operation 1, the resulting plane rotations are denoted by $P^{[c]}$ and $\dot{P}^{[c]}$.
We define $({\widetilde{H}}^{[0]},{\widetilde{K}}^{[0]})=(P_{s_\sigma}{\widetilde{H}},P_{s_\sigma}{\widetilde{K}})$, then by $({\widetilde{H}}^{[c]},{\widetilde{K}}^{[c]})=(P^{[c]}{\widetilde{H}}^{[c-1]}\dot{P}^{[c]},P^{[c]}{\widetilde{K}}^{[c-1]}\dot{P}^{[c]})$, we denote the pencil which is unitarily equivalent to $({\widetilde{H}}^{[c-1]},{\widetilde{K}}^{[c-1]})$, where the $c$-th column of ${\widetilde{H}}^{[c-1]}$ and ${\widetilde{K}}^{[c-1]}$ is restored to Hessenberg structure. Thus, the product of these matrices $\prod_{c=1}^{m-2}P^{[m-c-1]}$ and $\prod_{c=1}^{m-2}\dot{P}^{[c]}$ form the unitary matrices which restore the Hessenberg structure of $({\widetilde{H}},{\widetilde{K}})$. The unitary matrices $P^{[c]}$ and $\dot{P}^{[c]}$ are the product of plane rotations.

The following example illustrates our column by column approach.
\begin{example}\label{ex2}
	Continuing Example \ref{ex1} from \eqref{eq:struct_Ex3.1}, the structure of ${\widetilde{H}}^{[c]}$ and ${\widetilde{K}}^{[c]}$ is shown for $c=1,\dots,4$. In each step the entries marked $\textcolor{red}{\star}$ denote the entries which are used to construct the plane rotations establishing $P^{[c]} \in \mathbb{C}^{m \times m}$, the matrix
	which eliminates the entries below the first subdiagonal of the $(c+1)$-th column of ${\widetilde{H}}^{[c]}$ and ${\widetilde{K}}^{[c]}$ by unitary equivalence transformation, i.e., $({\widetilde{H}}^{[c+1]},{\widetilde{K}}^{[c+1]})=(P^{[c+1]}{\widetilde{H}}^{[c]}\dot{P}^{[c+1]},P^{[c+1]}{\widetilde{K}}^{[c]}\dot{P}^{[c+1]})$ has Hessenberg structure in its first $(c+1)$ columns.
Considering $({\widetilde{H}}^{[0]},{\widetilde{K}}^{[0]})=(P_{1}{\widetilde{H}},P_{1}{\widetilde{K}})$, we work column by column, starting from the first column up to the $4$-th column. Following the procedure outlined above, for the first three columns, the entries of the rows $r=5,6$ should be eliminated.
The structures that occur while restoring the first column to Hessenberg structure are
	\begin{eqnarray*}
	\scriptscriptstyle\substack{\begin{bmatrix}
\times & \times & \times & \times & \times &   \\
\textcolor{red}{\star} & \times & \times & \times &   &   \\
  & \times & \times & \times &   &   \\
  &   & \times & \times &   &   \\
\textcolor{red}{\star} & \times & \times & \times & \textcolor{red}{\star} &   \\
  &   &   &   & \times & \times
\end{bmatrix} \\{\widetilde{H}}^{[0]}}\rightarrow
\begin{bmatrix}
		\times & \times & \times & \times & \times &   \\
\textcolor{red}{\star} & \times & \times & \times & \times &   \\
  & \times & \times & \times &   &   \\
  &   & \times & \times &   &   \\
  & \times & \times & \times & \times &   \\
\textcolor{red}{\star} &   &   &   & \times &\textcolor{red}{\star}
	\end{bmatrix}\rightarrow \substack{\begin{bmatrix}
\times & \times & \times & \times & \times & \times \\
\times & \times & \times & \times & \times & \times \\
  & \times & \times & \times &   &   \\
  &   & \times & \times &   &   \\
  & \times & \times & \times & \times &   \\
  & \times & \times & \times & \times & \times
\end{bmatrix} \\{\widetilde{H}}^{[1]}},
\end{eqnarray*}
\begin{eqnarray*}
\scriptscriptstyle \substack{\begin{bmatrix}
\times & \times & \times & \times & \times &   \\
\textcolor{red}{\star} & \times & \times & \times &   &   \\
  & \times & \times & \times &   &   \\
  &   & \times & \times &   &   \\
\textcolor{red}{\star} & \times & \times & \times & \textcolor{red}{\star} &   \\
  &   &   &   &   & \times
\end{bmatrix} \\{\widetilde{K}}^{[0]} }\rightarrow
\begin{bmatrix}
	\times & \times & \times & \times & \times &   \\
\textcolor{red}{\star}& \times & \times & \times & \times &   \\
  & \times & \times & \times &   &   \\
  &   & \times & \times &   &   \\
  & \times & \times & \times & \times &   \\
  &   &   &   &   & \textcolor{red}{\star}
\end{bmatrix}\rightarrow \substack{\begin{bmatrix}
\times & \times & \times & \times & \times & \times \\
\times & \times & \times & \times & \times & \times \\
  & \times & \times & \times &   &   \\
  &   & \times & \times &   &   \\
  & \times & \times & \times & \times &   \\
  & \times & \times & \times & \times & \times
\end{bmatrix}\\ {\widetilde{K}}^{[1]}},
\end{eqnarray*}
where $({\widetilde{H}}^{[1]},{\widetilde{K}}^{[1]})=(P^{[1]}{\widetilde{H}}^{[0]}\dot{P}^{[1]},P^{[1]}{\widetilde{K}}^{[0]}\dot{P}^{[1]})$, and $P^{[1]}$ and $\dot{P}^{[1]}$ are the product of two plane rotations obtained by Operation 1. For the second column the intermediate structures are
	\begin{eqnarray*}
	\scriptscriptstyle \substack{\begin{bmatrix}
\times & \times & \times & \times & \times & \times \\
\times & \times & \times & \times & \times & \times \\
  & \textcolor{red}{\star} & \times & \times &   &   \\
  &   & \times & \times &   &   \\
  & \textcolor{red}{\star} & \times & \times & \textcolor{red}{\star} &   \\
  & \times & \times & \times & \times & \times
\end{bmatrix} \\{\widetilde{H}}^{[1]}}\rightarrow
	\begin{bmatrix}
\times & \times & \times & \times & \times & \times \\
\times & \times & \times & \times & \times & \times \\
  & \textcolor{red}{\star} & \times & \times & \times &   \\
  &   & \times & \times &   &   \\
  &   & \times & \times & \times &   \\
  & \textcolor{red}{\star} & \times & \times & \times & \textcolor{red}{\star}
\end{bmatrix}\rightarrow  \substack{\begin{bmatrix}
\times & \times & \times & \times & \times & \times \\
\times & \times & \times & \times & \times & \times \\
  & \times & \times & \times & \times & \times \\
  &   & \times & \times &   &   \\
  &   & \times & \times & \times &   \\
  &   & \times & \times & \times & \times
\end{bmatrix}\\{\widetilde{H}}^{[2]}},
\end{eqnarray*}
\begin{eqnarray*}
	\scriptscriptstyle \substack{\begin{bmatrix}
\times & \times & \times & \times & \times & \times \\
\times & \times & \times & \times & \times & \times \\
  & \textcolor{red}{\star} & \times & \times &   &   \\
  &   & \times & \times &   &   \\
  & \textcolor{red}{\star} & \times & \times & \textcolor{red}{\star} &   \\
  & \times & \times & \times & \times & \times
\end{bmatrix}\\ {\widetilde{K}}^{[1]}}\rightarrow
	\begin{bmatrix}
\times & \times & \times & \times & \times & \times \\
\times & \times & \times & \times & \times & \times \\
  & \textcolor{red}{\star} & \times & \times & \times &   \\
  &   & \times & \times &   &   \\
  &   & \times & \times & \times &   \\
  & \textcolor{red}{\star} & \times & \times & \times & \textcolor{red}{\star}
\end{bmatrix}\rightarrow \substack{\begin{bmatrix}
\times & \times & \times & \times & \times & \times \\
\times & \times & \times & \times & \times & \times \\
  & \times & \times & \times & \times & \times \\
  &   & \times & \times &   &   \\
  &   & \times & \times & \times &   \\
  &   & \times & \times & \times & \times
\end{bmatrix} \\ {\widetilde{K}}^{[2]}}.
\end{eqnarray*}
Now, we proceed by eliminating the entries below the subdiagonal of the third column.
	\begin{eqnarray*}
	\scriptscriptstyle \substack{\begin{bmatrix}
\times & \times & \times & \times & \times & \times \\
\times & \times & \times & \times & \times & \times \\
  & \times & \times & \times & \times & \times \\
  &   & \textcolor{red}{\star}   & \times &   &   \\
  &   & \textcolor{red}{\star}   & \times & \textcolor{red}{\star}   &   \\
  &   & \times & \times & \times & \times
\end{bmatrix}\\ {\widetilde{H}}^{[2]}}\rightarrow
	\begin{bmatrix}
\times & \times & \times & \times & \times & \times \\
\times & \times & \times & \times & \times & \times \\
  & \times & \times & \times & \times & \times \\
  &   & \textcolor{red}{\star} & \times & \times &   \\
  &   &   & \times & \times &   \\
  &   & \textcolor{red}{\star}& \times & \times & \textcolor{red}{\star}
\end{bmatrix}\rightarrow \substack{\begin{bmatrix}
\times & \times & \times & \times & \times & \times \\
\times & \times & \times & \times & \times & \times \\
  & \times & \times & \times & \times & \times \\
  &   & \times & \times & \times & \times \\
  &   &   & \times & \times &   \\
  &   &   & \times & \times & \times
\end{bmatrix} \\ {\widetilde{H}}^{[3]}},
\end{eqnarray*}
\begin{eqnarray*}
	\scriptscriptstyle \substack{\begin{bmatrix}
\times & \times & \times & \times & \times & \times \\
\times & \times & \times & \times & \times & \times \\
  & \times & \times & \times & \times & \times \\
  &   & \textcolor{red}{\star} & \times &   &   \\
  &   & \textcolor{red}{\star}  & \times & \textcolor{red}{\star}  &   \\
  &   & \times & \times & \times & \times
\end{bmatrix} \\ {\widetilde{K}}^{[2]}}\rightarrow
 \begin{bmatrix}
\times & \times & \times & \times & \times & \times \\
\times & \times & \times & \times & \times & \times \\
  & \times & \times & \times & \times & \times \\
  &   & \textcolor{red}{\star} & \times & \times &   \\
  &   &   & \times & \times &   \\
  &   & \textcolor{red}{\star} & \times & \times & \textcolor{red}{\star}
\end{bmatrix}\rightarrow \substack{\begin{bmatrix}
\times & \times & \times & \times & \times & \times \\
\times & \times & \times & \times & \times & \times \\
  & \times & \times & \times & \times & \times \\
  &   & \times & \times & \times & \times \\
  &   &   & \times & \times &   \\
  &   &   & \times & \times & \times
\end{bmatrix} \\ {\widetilde{K}}^{[3]}}.
\end{eqnarray*}
Finally, the entry at $r=6$, $c=4$ is eliminated and we obtain a Hessenberg pencil
	\begin{eqnarray*}
	\scriptscriptstyle \substack{\begin{bmatrix}
\times & \times & \times & \times & \times & \times \\
\times & \times & \times & \times & \times & \times \\
  & \times & \times & \times & \times & \times \\
  &   & \times & \times & \times & \times \\
  &   &   &  \textcolor{red}{\star}& \times &   \\
  &   &   &  \textcolor{red}{\star} & \times &  \textcolor{red}{\star}
\end{bmatrix} \\ {\widetilde{H}}^{[3]}}\rightarrow
\substack{\begin{bmatrix}
\times & \times & \times & \times & \times & \times \\
\times & \times & \times & \times & \times & \times \\
  & \times & \times & \times & \times & \times \\
  &   & \times & \times & \times & \times \\
  &   &   & \times & \times & \times \\
  &   &   &   & \times & \times
\end{bmatrix} \\ {\widetilde{H}}^{[4]}},
\end{eqnarray*}
\begin{eqnarray*}
	\scriptscriptstyle \substack{\begin{bmatrix}
\times & \times & \times & \times & \times & \times \\
\times & \times & \times & \times & \times & \times \\
  & \times & \times & \times & \times & \times \\
  &   & \times & \times & \times & \times \\
  &   &   & \textcolor{red}{\star} & \times &   \\
  &   &   & \textcolor{red}{\star} & \times & \textcolor{red}{\star}
\end{bmatrix}  \\ {\widetilde{K}}^{[3]}}
\rightarrow \substack{\begin{bmatrix}
\times & \times & \times & \times & \times & \times \\
\times & \times & \times & \times & \times & \times \\
  & \times & \times & \times & \times & \times \\
  &   & \times & \times & \times & \times \\
  &   &   & \times & \times & \times \\
  &   &   &   & \times & \times
\end{bmatrix} \\ {\widetilde{K}}^{[4]}}.
\end{eqnarray*}
\end{example}
We now describe the general process illustrated in the last example, i.e., how to construct matrices $P^{[c]}$ and $\dot{P}^{[c]}$ which are, in fact, the products of plane rotations. Without loss of generality, we assume that $(m-1)-(s_{\sigma}+1)$ poles have been added in the previous steps, and in the $\sigma$-th step, we restore the Hessenberg pencil such that these poles are conserved. For this purpose, we exploit Operation 1 to eliminate $((m-1)-(s_{\sigma}+1))\cdot (s_{\sigma}+1)+\sum_{j=1}^{\sigma}s_{j}$ entries from $P_{s_\sigma}\widetilde{H}$ and $P_{s_\sigma}\widetilde{K}$ to restore the Hessenberg structure, and for each element located on position $(r,c)$ two rotations $P$ and $\dot{P}$ are required. We require these elements $h_{c+1,c},h_{c+1,r},h_{r,c},h_{r,r}$ and $k_{c+1,c},k_{c+1,r},k_{r,c},k_{r,r}$ of the last computed pencil to be stored in $2\times 2$ blocks $A$ and $B$, introduced in Operation 1, to create these two rotations.

Following the elimination of all entries below the first subdiagonal, the matrix $\widetilde{J}$ can be represented by the following (restored) Hessenberg pencil
\begin{equation*}
	\widetilde{J}\widetilde{Q}P^{H}_{s_\sigma}\Big(\prod_{c=1}^{m-2}P^{[m-c-1]}\Big)^{H} {\widetilde{K}}^{[m-2]}=\widetilde{Q}P^{H}_{s_\sigma}\Big(\prod_{c=1}^{m-2}P^{[m-c-1]}\Big)^{H} {\widetilde{H}}^{[m-2]}.
\end{equation*}

Thus, the required structure, a Hessenberg pencil, is now
obtained by these unitary equivalence transformations. However, the ratio of subdiagonal elements of the last $s_{\sigma}+1$ rows will, in general, not be equal to the new poles that should be introduced in step $\sigma$. In the next section, we describe how these poles can be added to the subdiagonal of the Hessenberg pencil.

\subsection{Adding and swapping poles}\label{sec:addSwapPoles}
If we consider $({\widetilde{H}}^{[m-2]},{\widetilde{K}}^{[m-2]})$ as the Hessenberg pencil at the step $\sigma$, we need to add $s_{\sigma}+1$ poles to the entries of the pencil on position $(r,c)$ with
\begin{equation*}
	(r,c)=(m-s_{\sigma},m-s_{\sigma}-1),(m-s_{\sigma}+1,m-s_{\sigma}),\ldots,(m,m-1).
\end{equation*}
The first pole can be added to the last position, i.e., $(m,m-1)$. Let us consider $h_{m,m-1},h_{m,m}$ and $k_{m,m-1},k_{m,m}$ as the elements on the positions $(m,m-1)$ and $(m,m)$, respectively, in the matrices ${\widetilde{H}}^{[m-2]}$ and ${\widetilde{K}}^{[m-2]}$. Since $h_{m,m}$ and $k_{m,m}$ are both non-zero, and the matrices ${\widetilde{H}}^{[m-2]}$ and ${\widetilde{K}}^{[m-2]}$ form an unreduced Hessenberg pencil, the first pole $\frac{\mu_{m-s_{\sigma}-1}}{\nu_{m-s_{\sigma}-1}}=\psi_{m-s_{\sigma}-1}$ can be added through the following operation.
\begin{itemize}
    \item{\textbf{Operation 2} (Adding a pole).}

To add a pole to an unreduced Hessenberg pencil \cite{CaThesis,CaMeVa19_QZ}, let us consider the non-zero elements $h_{m,m-1},h_{m,m}$ and $k_{m,m-1},k_{m,m}$ on the positions $(m,m-1)$ and $(m,m)$ respectively, in the Hessenberg matrices ${\widetilde{H}}^{[m-2]}$, ${\widetilde{K}}^{[m-2]} \in \mathbb{C}^{m \times m}$. A plane rotation $\dot{P}=\begin{bmatrix}\bar{c}&-\bar{d}\\
	d&c\end{bmatrix}$ exists such that
\begin{equation*}
	\begin{bmatrix}
		h_{m,m-1}&h_{m,m}
	\end{bmatrix}\begin{bmatrix}
		\bar{c}\\
		d
	\end{bmatrix}=\mu_{m-s_{\sigma}-1},\quad \quad
	\begin{bmatrix}
		k_{m,m-1}&k_{m,m}
	\end{bmatrix}\begin{bmatrix}
		\bar{c}\\
		d
	\end{bmatrix}=\nu_{m-s_{\sigma}-1},
\end{equation*}
with $\frac{\mu_{m-s_{\sigma}-1}}{\nu_{m-s_{\sigma}-1}}=\psi_{m-s_{\sigma}-1}$.
\end{itemize}

Then the pole is swapped to the correct position $(m-s_{\sigma},m-s_{\sigma}-1)$ via the pole swapping operation.

\begin{itemize}
    \item{\textbf{Operation 3} (Swapping two consecutive poles).}
Two consecutive poles in an unreduced Hessenberg pair can be swapped via a unitary equivalence transformation \cite{VD81,CaMeVa19_QZ,CaThesis}.
Consider a Hessenberg pencil $({\widetilde{H}},{\widetilde{K}}) \in \mathbb{C}^{m \times m}\times \mathbb{C}^{m \times m}$, where we will swap the poles $\psi = \frac{\tau}{\kappa}$ and $\xi = \frac{\mu}{\nu}$
\begin{eqnarray*}
	\left(\left[\begin{array}{cccccc}
		\times & \times & \cdots &\times &\times&  \times  \\
		\times & \times & \cdots  &\times&\times&     \times  \\
		& \times & \cdots   &\times&\times&  \times \\
		
		& &\ddots  &\vdots &\vdots&  \vdots                  	\\
		& &   &\tau&\ast&         \times        \\
		&  &   &&\mu&\times\\
	\end{array}\right], \left[\begin{array}{cccccc}
		 \times & \times & \cdots &\times &\times&\times\\
		\times	& \times &  \cdots  &\times &\times&\times\\
		& \times & \cdots &\times&\times&\times\\
		
		& &\ddots &\vdots&  \vdots&\vdots\\
			&  &   &\kappa&\ast&\times\\
		& &   &&\nu&\times\\
	\end{array}\right]\right).
\end{eqnarray*}
Suppose that $A$ and $B$ refer to the small $2 \times 2$ upper triangular matrices taken out of the large Hessenberg pair. We are thus looking for a unitary equivalence transformation of the following form
\begin{eqnarray}\label{eq223}
	P \Big(\begin{bmatrix}
		\tau & \ast\\
		0 & \mu
	\end{bmatrix} -\xi\begin{bmatrix}
		\kappa& \ast\\
		0 & \nu
	\end{bmatrix}\Big)\dot{P}=\begin{bmatrix}
		\hat{\tau} & \star\\
		0 & \hat{\mu}
	\end{bmatrix} -\xi\begin{bmatrix}
		\hat{\kappa}& \star\\
		0 & \hat{\nu}
	\end{bmatrix},
\end{eqnarray}
with $\frac{\hat{\mu}}{\hat{\nu}}=\frac{\tau}{\kappa}=\psi$, and $\frac{\hat{\tau}}{\hat{\kappa}}=\frac{\mu}{\nu}=\xi$. In this regard, $\dot{P}$ is first constructed in a way that $Z \dot{P}e_{1}=\begin{bmatrix} 0& 0\end{bmatrix}^{\top}$ with
\begin{equation*}
	Z=\nu A-\mu B=\begin{bmatrix}
		\times & \times\\
		0 & 0
	\end{bmatrix}.
\end{equation*}
This implies that $A\dot{P}e_{1}$ and $B\dot{P}e_{1}$ are parallel vectors and that a rotation $P$ can be computed to simultaneously introduce a zero in position $(2,1)$ of both $AP$ and $BP$, thus obtaining \eqref{eq223}.
\end{itemize}

Operation 3 can be used to swap any two poles in an unreduced Hessenberg pencil.
Thus, we can introduce a pole at any position in the pencil by first using Operation 2 to introduce it in the pencil and then apply Operation 3 to swap this pole to the correct position.

Let us consider the unitary matrices $P^{(i)}$ and $\dot{P}^{(i)}$ as the products of the plane rotations $P$ and $\dot{P}$ to add the $i$-th pole to $({\widetilde{H}}^{[m-2+(i-1)]},{\widetilde{K}}^{[m-2+(i-1)]})$ and swap it to the correct position $(m-s_{\sigma}-1+i,m-s_{\sigma}-2+i)$, which results in
$$({\widetilde{H}}^{[m-2+i]},{\widetilde{K}}^{[m-2+i]})=(P^{(i)}{\widetilde{H}}^{[m-2+(i-1)]}\dot{P}^{(i)},P^{(i)}{\widetilde{K}}^{[m-2+(i-1)]}\dot{P}^{(i)}).$$

In the following example we illustrate the process of adding and swapping poles in a Hessenberg pencil.
\begin{example}\label{ex3}
	We continue with the Hessenberg pencil $({\widetilde{H}}^{[4]},{\widetilde{K}}^{[4]})$ obtained in Example \ref{ex2}. Without loss of generality, we assume that three poles have been added to $({\hat{H}},{\hat{K}}) \in \mathbb{C}^{4 \times4} \times \mathbb{C}^{4 \times4}$ in previous steps, and as $s_{3}=1$, we require to add two poles at this stage. The ratio of subdiagonal elements marked $\textcolor{blue}{\star}$, $\textcolor{red}{\star}$ of ${\widetilde{H}}^{[6]}$ and ${\widetilde{K}}^{[6]}$ (specified subsequently) denote $\psi_{4}$ and $\psi_{5}$, respectively. In the following, we display the changes on the Hessenberg pencil $({\widetilde{H}}^{[4]},{\widetilde{K}}^{[4]})$.
\end{example}
We start by adding the pole $\psi_4$ to the pencil in position $(6,5)$ and swap it to the position $(5,4)$
\begin{eqnarray*}
	\scriptscriptstyle \substack{\begin{bmatrix}
			\times & \times & \times &\times &\times&\times\\
			\times& \times & \times  &\times&\times&\times\\
			& \times & \times   &\times&\times&\times\\
			&&\times  &\times &\times&\times\\
			& &  &\times&\times&\times\\
			&  &   &&\times&\times\\
		\end{bmatrix} \\ {\widetilde{H}}^{[4]}}\rightarrow
	\begin{bmatrix}
		\times & \times & \times &\times &\times&\times\\
		\times& \times & \times  &\times&\times&\times\\
		& \times & \times   &\times&\times&\times\\
		&&\times  &\times &\times&\times\\
		& &  &\times&\times&\times\\
		&  &   &&\textcolor{blue}{\star}&\times\\
	\end{bmatrix} \rightarrow \substack{\begin{bmatrix}
	\times & \times & \times &\times &\times&\times\\
	\times& \times & \times  &\times&\times&\times\\
	& \times & \times   &\times&\times&\times\\
	&&\times  &\times &\times&\times\\
	& &  &\textcolor{blue}{\star}&\times&\times\\
	&  &   &&\times&\times\\
\end{bmatrix}  \\ {\widetilde{H}}^{[5]}},
\end{eqnarray*}
\begin{eqnarray*}
	\scriptscriptstyle \substack{\begin{bmatrix}
			\times& \times & \times &\times &\times&\times\\
			\times & \times & \times &\times&\times&\times\\
			& \times & \times   &\times&\times&\times\\
			& &\times  &\times &\times&\times\\
			& &  & \times&\times&\times\\
			&  &  &&\times&\times\\
		\end{bmatrix} \\ {\widetilde{K}}^{[4]}}\rightarrow
	\begin{bmatrix}
		\times& \times & \times &\times &\times&\times\\
		\times & \times &\times  &\times&\times&\times\\
		& \times & \times   &\times&\times&\times\\
		& &\times  &\times &\times&\times\\
		& &  & \times&\times&\times\\
		&  &  &&\textcolor{blue}{\star}&\times\\
	\end{bmatrix}\rightarrow \substack{\begin{bmatrix}
	\times& \times & \times &\times &\times&\times\\
	\times & \times & \times &\times&\times&\times\\
	& \times & \times   &\times&\times&\times\\
	& &\times  &\times &\times&\times\\
	& &  &\textcolor{blue}{\star}&\times&\times\\
	&  &  &&\times&\times\\
\end{bmatrix}\\ {\widetilde{K}}^{[5]}},
\end{eqnarray*}
where $({\widetilde{H}}^{[5]},{\widetilde{K}}^{[5]})=(P^{(1)}{\widetilde{H}}^{[4]}\dot{P}^{(1)},P^{(1)}{\widetilde{K}}^{(4)}\dot{P}^{(1)})$, and $P^{(1)}$ and $\dot{P}^{(1)}$ are the product of two plane rotations. Eventually the last pole should be added
\begin{eqnarray*}
	\scriptscriptstyle \substack{\begin{bmatrix}
			\times & \times & \times &\times &\times&\times\\
			\times& \times & \times  &\times&\times&\times\\
			& \times & \times   &\times&\times&\times\\
			&&\times &\times &\times&\times\\
			& &  &\textcolor{blue}{\star}&\times&\times\\
			&  &   &&\times&\times\\
		\end{bmatrix} \\ {\widetilde{H}}^{[5]}}\rightarrow
\substack{\begin{bmatrix}
			\times & \times & \times &\times &\times&\times\\
			\times& \times & \times  &\times&\times&\times\\
			& \times & \times   &\times&\times&\times\\
			&&\times &\times &\times&\times\\
			& &  &\textcolor{blue}{\star}&\times&\times\\
			&  &   &&\textcolor{red}{\star}&\times\\
		\end{bmatrix} \\ {\widetilde{H}}^{[6]}},
\end{eqnarray*}
\begin{eqnarray*}
	\scriptscriptstyle  \substack{\begin{bmatrix}
			\times& \times & \times &\times &\times&\times\\
			\times &\times  & \times &\times&\times&\times\\
			& \times & \times   &\times&\times&\times\\
			& &\times  &\times &\times&\times\\
			& &  &\textcolor{blue}{\star}&\times&\times\\
			&  &  &&\times&\times\\
		\end{bmatrix}\\ {\widetilde{K}}^{[5]}}\rightarrow
	 \substack{ \begin{bmatrix}
			\times& \times & \times &\times &\times&\times\\
			\times & \times &\times  &\times&\times&\times\\
			& \times & \times   &\times&\times&\times\\
			& &\times  &\times &\times&\times\\
			& &  &\textcolor{blue}{\star}&\times&\times\\
			&  &  &&\textcolor{red}{\star}&\times\\
		\end{bmatrix} \\ {\widetilde{K}}^{[6]}}.
\end{eqnarray*}
In Examples \ref{ex1}, \ref{ex2}, and \ref{ex3}, we consider the continuous Gegenbauer-Sobolev inner product \eqref{eq:inprod_example}, and we discretize such that we can generate a sequence of two SORFs $\{r_{0},r_{1}\}$ with the pole $\psi_{1}=\xi_{1}$, which is already added to the Hessenberg pencil computed in the first step (i.e., corresponding to $s_{1}=1$ with the Jordan matrix of size $2\times2$). The remaining poles $\{\psi_{k}\}_{k=2}^{5}$ can be freely chosen.

The process of swapping any two consecutive poles in an unreduced Hessenberg pair via a unitary equivalence on itself is done via Operation 3. In general, after adding the first pole, we repeat the process of swapping to transfer the first pole to the position $(m-s_{\sigma},m-s_{\sigma}-1)$. The number of swapping for the first pole of step $\sigma$ is $s_{\sigma}$, and generally for the pole $\psi_{k}$, $k=1,2 \ldots,s_{\sigma}+1$, the process of swapping should be repeated $s_{\sigma}-k+1$ times. To swap two arbitrary consecutive poles $\xi=\frac{h_{r,c}}{k_{r,c}}$ and $\psi=\frac{h_{r-1,c-1}}{k_{r-1,c-1}}$ in the last computed Hessenberg pencil, these elements $h_{r-1,c-1},h_{r-1,c},h_{r,c-1},h_{r,c}$ and $k_{r-1,c-1},k_{r-1,c},k_{r,c-1},k_{r,c}$ are required to be stored in $2 \times 2$ blocks $A$ and $B$, given in Operation 3.

We finally get the following relation
\begin{multline*}
	\widetilde{J}\widetilde{Q}P^{H}_{s_\sigma}\Big(\prod_{c=1}^{m-2}P^{[m-c-1]}\Big)^{H}\Big(\prod_{i=1}^{s_{\sigma}+1}P^{(s_{\sigma}-i+2)}\Big)^{H}{\widetilde{K}}^{[m+s_{\sigma}-1]}\\
	=\widetilde{Q}P^{H}_{s_\sigma}\Big(\prod_{c=1}^{m-2}P^{[m-c-1]}\Big)^{H}\Big(\prod_{i=1}^{s_{\sigma}+1}P^{(s_{\sigma}-i+2)}\Big)^{H}{\widetilde{H}}^{[m+s_{\sigma}-1]},
\end{multline*}
thereby by defining
\begin{equation*}
\widetilde{Q}=\widetilde{Q}P^{H}_{s_\sigma}\Big(\prod_{c=1}^{m-2}P^{[m-c-1]}\Big)^{H}\Big(\prod_{i=1}^{s_{\sigma}+1}P^{(s_{\sigma}-i+2)}\Big)^{H}, \quad {\widetilde{H}}={\widetilde{H}}^{[m+s_{\sigma}-1]},\quad {\widetilde{K}}={\widetilde{K}}^{[m+s_{\sigma}-1]},
\end{equation*}
the Hessenberg pencil IEP \ref{IEP} is solved.

\subsection{Alternative strategy}\label{sec:SORFviaSOP}
Above, we introduced an algorithm that solves Problem \ref{IEP} by splitting the discrete inner product into blocks and updating a known solution when a single block is added. The poles are added to each block, thus after adding a block we have a Hessenberg pencil which encodes a sequence of rational functions.

Alternatively, we can first generate the full sequence of Sobolev orthogonal polynomials with respect to the discretized inner product \eqref{eq:DiscretizedInProd} using the methods in \cite{VB23,GaZh95}, this leads to a Hessenberg matrix $\widetilde{H}$ of size $m\times m$.
Second, we notice that this Hessenberg matrix can be interpreted as the Hessenberg pencil $(\widetilde{H},\widetilde{K})$, with $\widetilde{K} = I_m$ the $m\times m$ identity matrix. This pencil implies that all poles are infinity (since the subdiagonal of $I_m$ contains all zero entries), and therefore contains the recurrence coefficients for a sequence of SOPs orthogonal with respect to \eqref{eq:DiscretizedInProd}.
Third, using Operation 2 and Operation 3 from Section \ref{sec:addSwapPoles}, we can introduce the poles in $\Xi = \{\xi_k\}_{k=1}^{N-1}$ to the leading principal $N\times N$ subpencil. This $N\times N$ subpencil then contains the recurrence coefficients for generating the SORFs orthogonal with respect to \eqref{eq:SobInprod}.

\section{Numerical experiments}	\label{sec:numerics}
In this section we assess the proposed updating procedure and the procedure described in Section \ref{sec:SORFviaSOP} by generating the SORFs numerically, and compare them to the Krylov subspace method \cite{FaVBVBVa24,Ru94}.

Consider the Jordan-like matrix $J \in \mathbb{C}^{m \times m}$, constructed by distinct nodes $\{z_{j}\}_{j=1}^{\sigma}$ and derivative orders $\{s_j\}_{j=1}^\sigma$, the weight vector $w \in \mathbb{C}^{m}$, and the set of poles $\Xi = \{\xi_{k}\}_{k=1}^{N-1}$, and $\Psi = \{\psi_{k}\}_{k=1}^{m-1}$ with $\psi_{k}=\xi_{k}$, $k=1,2,\ldots,N-1$, given in Problem \ref{IEP}. The solution strategies, the two procedures introduced in this paper and the Krylov subspace method \cite{FaVBVBVa24}, compute a solution in the form of a Hessenberg pencil $(H_{m},~K_{m}) \in\mathbb{C}^{m\times m}$, with subdiagonal entries satisfying
$\frac{h_{k+1,k}}{k_{k+1}}= \psi_{k} $, $k=1,2,\ldots,m-1$, and $m=\sum_{j=1}^\sigma (s_{j}+1)$ such that
\[
	Q_{m}^{H} J Q_{m} K_{m}=H_{m},
	\]
	where $Q_{m}^{H} Q_{m} = I_{m}$ and $Q_{m} e_{1} = w/|| w ||_{2}$. The following metrics for the errors of the computed pencil are used.
\begin{itemize}
\item The accuracy of the recurrence relation, consisting of the recurrence
matrices $(H_{m},~K_{m})$ and basis $Q_{m}$, is measured by
\begin{equation*}
E_{m}(r)=\dfrac{\big\Vert J Q_{m} K_{m}-Q_{m}H_{m}\big\Vert_{2}}{\max \big(\big\Vert J Q_{m} K_{m}\big\Vert_{2},\big\Vert Q_{m} H_{m}\big\Vert_{2}\big)}.
\end{equation*}

\item The following error metric quantifies the accuracy of the poles. The poles of the
computed pencil $(H_{m},~K_{m})$ are compared to the given poles $\{\psi_{k}\}_{k=1}^{m-1}$,
\begin{equation*}
	E_{m}(p)=\max_{k=1,2,\ldots,m-1}\left\{\frac{\left|\frac{H_{m}(k+1,k)}{K_{m}(k+1,k)}-\psi_{k}\right|}{\vert \psi_{k}\vert} \right\}.
\end{equation*}

\item The orthonormality of the formed basis $Q_{m}$ is measured by
	\begin{equation*}
	E_{m}(Q) = \big\Vert Q_{m}^{H} Q_{m} - I_{m}\big\Vert_{2}.
	\end{equation*}

\item As a metric for orthonormality of the rational functions, we compute the two-norm of the difference between the moment matrix and the identity matrix
\begin{equation*}
    E_{m}^x(S) = \Vert M^x_m - I_m\Vert_2.
\end{equation*}
For the discretized Sobolev inner product, the moment matrix $M^d_m = \{\langle r_i,r_j\rangle_S\}_{i,j=0}^{m-1}$ requires evaluations of the rational functions (and its derivatives) at the quadrature nodes.
A sequence of SORFs $\{r_{k} \}_{k=0}^{m-1}$ can be evaluated in the point $z$ by using its recurrence pencil $(H_{m},~K_{m})$ (see \cite{FaVBVBVa24} for details).
For the continuous Sobolev inner product, we use Clenshaw-Curtis quadrature to compute each element of the moment matrix $M^c_m = \{(r_i,r_j)_S\}_{i,j=0}^{m-1}$.

\end{itemize}


In Section \ref{sec:Num_exactness}, we illustrate the results obtained in Section \ref{sec:theory}, most importantly the relationship between the continuous and discretized inner product in terms of the degree of exactness.
Then, in Section \ref{sec:Num_comparison} a comparison between three numerical algorithms for the generation of SORFs is made.

\subsection{Degree of exactness}\label{sec:Num_exactness}
In order to generate the first $N$ SORFs with respect to a continuous Sobolev inner product, a discretized inner product with $\sigma = (s+1)(N-1) +1$ nodes is used here, and the corresponding $m\times m$ recurrence pencil is computed.
The recurrence coefficients for SORFs with respect to the continuous inner product are in the $N\times N$ leading principal submatrix of the recurrence pencil generated for the discretized inner product.
We illustrate this by the moment matrices associated with the continuous and discretized inner product.

\begin{figure}[H]
   \begin{subfigure}[b]{0.45\textwidth}
        \centering
		\setlength\figureheight{5cm}
		\setlength\figurewidth{12.5cm}
        \input{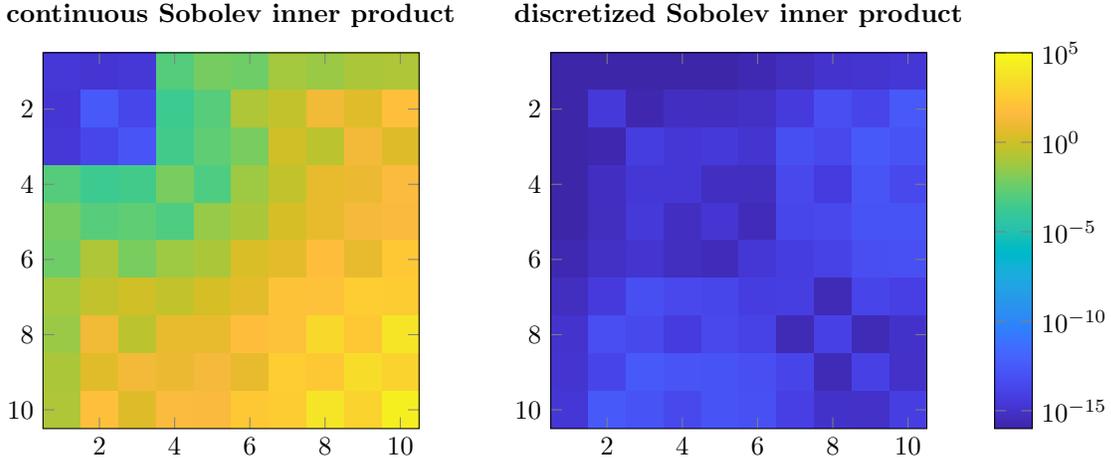}
    \end{subfigure}
    \hspace{0.05\textwidth}
    \caption{Orthogonality error for each entry of the moment matrix (minus the identity matrix of the same size) generated by the continuous (left) and discretized (right) Sobolev-Gegenbauer inner product. The first 3 SORFs with respect to the continuous inner product are computed, which require the computation of a $10\times 10$ Hessenberg pencil. As predicted by the degree of exactness of the discretization, the orthogonality error for the continuous inner product increases rapidly beyond the $3\times 3$ leading principal subpencil, whereas the orthogonality error for the discretized inner product is small everywhere.}\label{fig:momentMatrix_cont_vs_disc}

    \end{figure}

As the continuous inner product, we take the Gegenbauer-Sobolev inner product \eqref{eq:inprod_example} and we discretize such that we can generate the first 3 SORFs, $\{r_0,r_1,r_2\}$ with poles $\Xi = \{-\omega,\omega\}$ (Example \ref{ex:gegenbauer_poles} with $M=1$). Since $s = 1$ and $N=3$, we require $\sigma = (s+1)(N-1)+1 =  5$ quadrature nodes and from Example \ref{ex:gegenbauer_matrixJ} we know that $J$ is of size $10\times 10$ and $w$ of size $10 \times 1$, thus $m=10$.
Using the updating procedure, we generate the Hessenberg pencil $(H_m,K_m)$ and place the poles $-\omega$ on position $(2,1)$ and $\omega$ on position $(3,2)$.
We compute the moment matrix $M^c_{10}$ for the continuous inner product \eqref{eq:inprod_example} and $M^d_{10}$ for the discretized inner product \eqref{eq:discretizedGegenbauerSobolevInprod} and, in Figure \ref{fig:momentMatrix_cont_vs_disc}, we show the size of the entries of the matrices $M^c_{10}-I_{10}$ and $M^d_{10}-I_{10}$, which can be viewed as a metric for the orthogonality error.
The largest orthogonality error for the discretized inner product is of size $10^{-12}$, thus the algorithm computed the recurrence pencil accurately.
For the continuous inner product we notice that the orthogonality error inside the $3\times 3$ leading principal submatrix is of size $10^{-12}$. However, outside this submatrix the orthogonality error increases rapidly, this is not due to any computational error, but merely the consequence of the fact that the equality $(f,g)_S = \langle f,g\rangle_S$ only holds for $f,g\in\textrm{span}\{r_0,r_1,r_2, r_1',r_2'\}$.
This verifies the results in Section \ref{sec:theory} which state that in order to generate $N$ SORFs we must compute a larger Hessenberg pencil whose size is related to the discretized inner product.

\subsection{Comparison of three algorithms}\label{sec:Num_comparison}
We compare the performance of three algorithms for generating SORFs with respect to the Sobolev-Gegenbauer inner product \eqref{eq:inprod_example} with poles as in Example \ref{ex:gegenbauer_poles}.
These three algorithms are the Krylov subspace method in \cite{FaVBVBVa24}, the updating procedure proposed in Section \ref{sec:algorithm}, and the alternative procedure in Section \ref{sec:SORFviaSOP}.

We track the error metrics introduced above as we increase $N$, the number of SORFs $\{r_k\}_{k=0}^{N-1}$.
The first SORF in the sequence, $r_0$ is a constant and requires one quadrature node, which leads to a $2\times 2$ Hessenberg pencil IEP through the discretized inner product.
Thereafter, the addition of one SORF, $r_{\ell-1}$, to the sequence $\{r_k\}_{k=0}^{\ell-2}$ requires two additional quadrature nodes, i.e., the size of the matrices in the IEP is $(2+(\ell-1)4)\times (2+(\ell-1)4)$. Let $m = (2+(\ell-1)4)$.

\begin{figure}[H]
   \begin{subfigure}[b]{0.45\textwidth}
        \centering
		\setlength\figureheight{8cm}
		\setlength\figurewidth{13cm}
        \input{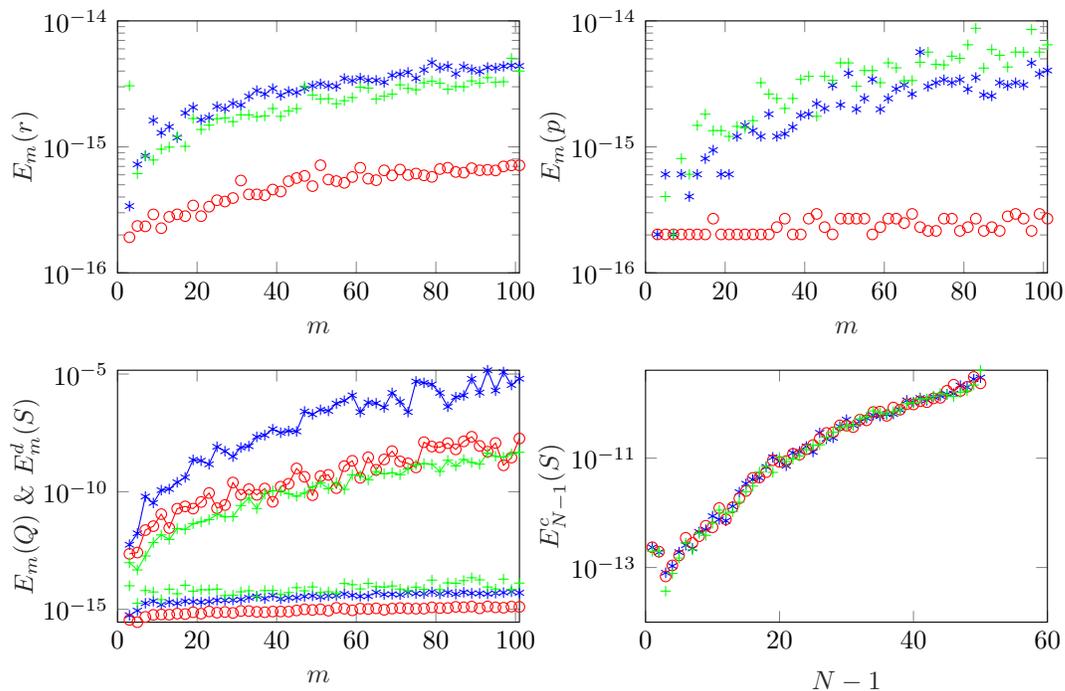}
    \end{subfigure}
    \hspace{0.05\textwidth}
    \caption{Error metrics for increasing the size of the sequence of SORFs requested. SORFs for the Gegenbauer-Sobolev inner product with $\lambda=1$ and $\mu=2$. The poles are $\{-\omega,\omega,-2\omega,2\omega,\dots\}$ with $\omega = 1.1$. The error metrics for the recurrence relation $E_m(r)$ and for the poles $E_m(p)$ indicate high accuracy. Krylov subspace method ({\color{red} $\circ$}), Updating procedure ({\color{blue} $\ast$}), and  the procedure in Section \ref{sec:SORFviaSOP} ({\color{green} +}).}\label{fig:momentMatrix_cont_vs_disc_1}
    \end{figure}

For the parameters $\mu = 2$, $\lambda=1$ we show the error metrics for $\omega = 1.1$ in Figure \ref{fig:momentMatrix_cont_vs_disc_1}.

The two top panels of Figure \ref{fig:momentMatrix_cont_vs_disc_1} show a recurrence error and pole error that is below $10^{-14}$ for all three algorithms, and the Krylov subspace method is the most accurate.

For the orthogonality, we track three metrics in the bottom two panels of the figure.
The following two metrics are shown in the lower left panel of the figure.
First, the orthonormality of the basis $Q_m\in\mathbb{R}^{m\times m}$ with respect to the Euclidean inner product, which is below $10^{-14}$ for all algorithms.
Second, we track $E_m^d(S)$ (data points connected by a full line), which measures the orthogonality of the full sequence of $m$ computed rational functions with respect to the discretized inner product \eqref{eq:discretizedGegenbauerSobolevInprod}. For this metric we observe that the Krylov subspace method and the approach from Section \ref{sec:SORFviaSOP} perform better than the updating procedure.
Although mathematically they should be the same, we observe a significant difference in practice. This can be partially explained by the computational error incurred when constructing the moment matrix $M_m^d$ for the discretized inner product by evaluating rational functions.

Third, we track the orthogonality of the first $N-1$ computed rational functions with respect to the continuous inner product \eqref{eq:SobInprod}. The corresponding metric $E_{N-1}^c(S)$ is shown in the lower right panel of the figure.
We observe that $E_{N-1}^c(S)$ remains below $10^{-11}$, which is smaller than the orthogonality error we observe for the corresponding $E_m^c(S)$. This difference is mainly due to the fact that the moment matrix $M_m^d$ is larger than the moment matrix $M_m^c$, and this larger size can lead to a larger accumulation of rounding error and larger errors in evaluating rational functions of higher degrees.

\section{Conclusion}
The problem of generating a sequence of rational functions orthogonal with respect to a Sobolev inner product is reformulated as a matrix problem, a Hessenberg pencil inverse eigenvalue problem.
We show how rational Gauss quadrature rules can be used to formulate this Hessenberg pencil inverse eigenvalue problem.
Solving this problem results in a Hessenberg pencil containing the recurrence coefficients of the sequence of rational functions.

Two algorithms are proposed to solve this inverse eigenvalue problem. The first is an updating procedure which is predominantly based on unitary similarity transformations performed on a matrix pencil.
The second starts by generating a sequence of polynomials orthogonal with respect to a discretized version of the Sobolev inner product, which results in a Hessenberg recurrence matrix. Using operations discussed in this paper, the Hessenberg matrix can be altered to a Hessenberg pencil that corresponds to a sequence of rational functions.
A numerical comparison between these two algorithms and a Krylov subspace method, introduced in another paper, is performed.

\section*{Funding}
The research was partially supported by the Research Council KU Leuven (Belgium), project C16/21/002 (Manifactor: Factor Analysis for Maps into Manifolds) and by the Fund for Scientific Research -- Flanders (Belgium), projects G0A9923N (Low rank tensor approximation techniques for up- and downdating of massive online time series clustering) and G0B0123N (Short recurrence relations for rational Krylov and orthogonal rational functions inspired by modified moments).
\section*{Declarations}
\textbf{Conflict of interest:} The authors declare that they have no conflict of interest.

\bibliographystyle{siam}
\bibliography{references}

\end{document}